\documentclass[12pt]{amsart}

\usepackage[utf8]{inputenc}   
\usepackage[T1]{fontenc}    
\usepackage[english]{babel}  
\usepackage{hyperref}
\usepackage{geometry} 
\usepackage{bm,amsmath,amssymb,amsthm,latexsym}

\newtheorem{theorem}{Theorem}
\newtheorem{proposition}[theorem]{Proposition}
\newtheorem{fact}[theorem]{Fact}
\newtheorem{cor}[theorem]{Corollary}
\newtheorem{lmm}[theorem]{Lemma}
\newtheorem*{claim}{Claim}

\theoremstyle{definition}
\newtheorem{definition}[theorem]{Definition}
\newtheorem{remark}[theorem]{Remark}
\newtheorem{conj}[theorem]{Conjecture}
\newtheorem*{expl}{Example}
\def\bsp{\begin{expl}}
\def\ebsp{\end{expl}}
\def\beh{\begin{claim}}
\def\ebeh{\end{claim}}
\def\defn{\begin{definition}\upshape}
\def\edefn{\end{definition}}
\def\satz{\begin{theorem}}
\def\esatz{\end{theorem}}
\def\tats{\begin{fact}}
\def\etats{\end{fact}}
\def\kor{\begin{cor}}
\def\ekor{\end{cor}}
\def\rema{\begin{remark}}
\def\erema{\end{remark}}
\def\lem{\begin{lmm}}
\def\elem{\end{lmm}}
\def\bem{\begin{remark}}
\def\ebem{\end{remark}}
\def\verm{\begin{conj}}
\def\everm{\end{conj}}
\def\bew{\par\noindent{\em Proof: }}

\def\satzli{\begin{proposition}}
\def\esatzli{\end{proposition}}

\def\gen#1{\langle #1\rangle} 

\def\BB{\mathfrak B}
\def\p{\varphi}

\begin{document}

\title{Approximate subgroups}          
\author{Jean-Cyrille Massicot and Frank O. Wagner}
\address{ENS de Rennes; Avenue Robert Schuman, Campus de Ker Lann, 35170 Bruz,
France}
\address{Universit\'e Lyon 1; CNRS; Institut Camille Jordan UMR 5208, 21
avenue
Claude Bernard, 69622 Villeurbanne-cedex, France}
\email{massicot@math.univ-lyon1.fr}
\email{wagner@math.univ-lyon1.fr}
\keywords{approximate subgroup; definability; definable amenability}

\begin{abstract}
Given a definably amenable approximate subgroup $A$ of a (local) group in some
first-order structure, there is a type-definable subgroup $H$ normalized by $A$
and contained in $A^4$ such that every definable superset of $H$ has positive
measure.
\end{abstract}
\thanks{Partially supported by ValCoMo (ANR-13-BS01-0006)}
\date{23 June 2014}
\subjclass[2010]{11B30; 20N99; 03C98; 20A15}

\maketitle

\section*{Introduction}
Let $G$ be a group and $K>0$ an
integer, a subset $A\subseteq G$ closed under inverse is a {\em
$K$-approximate subgroup}
if there is a finite subset $E\subseteq G$ with $|E|\le K$ such that
$A^2=\{ab:a,b\in A\}\subseteq EA$. Then $A^n\subseteq E^{n-1}A$.

\smallskip

Following work of Hrushovski \cite{H12} and many others, Breuillard, Green and
Tao \cite{BGT12} have classified finite approximate subgroups of local groups
(see \cite{vdD13} for an excellent survey). In particular, they show that there
is an approximate subgroup $A^*\subseteq A^4$ and an actual
$A^*$-invariant subgroup $H^*\subseteq A^*$ such that\begin{itemize}
\item finitely many left translates of $A^*$ cover $A$, and
\item $\gen{A^*}/H^*$ is nilpotent.
\end{itemize}

The result and its proof are inspired not only by Gleason's and Yamabe's
solution of Hilbert's $5$th problem \cite{Gl52,Ya53} and its extension to the
local context by Goldbring \cite{Go10}, but also by
Gromov's Theorem on groups with polynomial growth, and is indeed a way to
generalize this theorem. The three articles \cite{H12,BGT12,vdD13} provide some
applications to geometric group theory. 

\smallskip

The proof proceeds by considering a non-principal ultraproduct of a sequence of
finite counterexamples $(A_n:n<\omega)$ with $|A_n|\to\infty$, giving rise to
a pseudofinite counterexample $A$. Then an $A$-invariant subgroup $H\subseteq
A^4$ is constructed such that $\gen A/H$ is locally compact. From
Yamabe's theorem on the approximation of locally compact groups by Lie groups,
it follows that there are suitable $A^*$ and $H^*$ such that $\gen{A^*}/H^*$ is
a real Lie group; using pseudofiniteness, the final result is obtained.

\smallskip

The construction of the locally compact quotient $\gen A/H$ and the Lie model
$\gen{A^*}/H^*$ was first shown by Hrushovski \cite{H12} by model-theoretic
means
inspired by and reminiscent of stability theory. Breuillard, Green and Tao use
instead a (subsequent) theorem of Sanders \cite{S10} from finite combinatorics,
constructing successively the traces of the definable supersets of $H$ on the
various finite approximate groups $A_n$. Using the ultraproduct construction,
the pseudofinite counting measure and the translation-invariant ideal of
measure
zero sets, Hrushovski's theorem allows to recover Sanders' result at least
qualitatively.

\subsection*{Definability}

The topology of $\gen A/H$ was constructed analytically in~\cite{BGT12}, but
has a natural model-theoretic interpretation already given in \cite{H12}.
Recall that a subset of the ultraproduct is {\em definable} if it is the set
of realizations of some first-order formula (usually involving quantifiers); it
is {\em type-definable} if it is given as the intersection of a countable (say)
family of definable sets. For instance, the centralizer of a group element $g$
is defined by the formula $xg=gx$, and if $G$ is a group defined by a formula
$\p(x)$, its centre $Z(G)$ is defined by the formula $\p(x)\land\forall
y\,(\p(y)\to xy=yx)$. On the other hand, the group generated by an
element $g$, or the centre of a type-definable group, are in general not even
type-definable. 

\smallskip
If $H$ is a type-definable normal subgroup of $\gen A$, it has
{\em bounded index} if any definable superset of $H$ covers any definable subset
of $\gen A$ in finitely many translates. We can then endow the quotient $\gen
A/H$ with the \emph{logic topology} whose proper closed subsets are precisely
those subsets whose preimage in $\gen A$ is type-definable; this will turn it
into a locally compact 
topological group.

\smallskip

Of course, (type-)definability strongly depends on the language: if we
{\em expand} the structure, for instance by adding predicates for certain
subsets, there will be more definable sets. While the group $H$ constructed by
Hrushovski is naturally type-definable in the structure given,
it only becomes so in Sanders' Theorem (either in the ultraproduct or in a
suitable version using a bi-invariant measure instead of cardinality) after
such an expansion of language. Hrushosvki, on the other hand, assumes
the existence of a bi-invariant S1 ideal (which should be thought of as the
ideal of sets of measure zero) which in addition is automorphism invariant; in
order for an ideal to become automorphism invariant, one would generally also
have to expand the language. Such an expansion does not matter much if the
structure to start with is arbitrary, but should be avoided if the
initial structure has particular model-theoretic properties one
wants to preserve.

\smallskip

For example, Eleftheriou and Peterzil \cite{EP12} construct $H$
type-definably without expanding the language in the case when $A$ is definable
in an $o$-minimal expansion of an ordered group (such as the field of real
numbers with exponentiation), provided that $\gen A$ is
abelian. Pillay \cite{Pi14}, generalizing an argument in \cite{HP09},
generalizes this result if $A$ is definable in a theory without the
independence property, and is definably amenable (see below). In fact, in this
setting there is a unique minimal choice for $H$,
namely the unique minimal type-definable subgroup of bounded index,
$\gen A^{00}$.

\subsection*{Definable amenability}

We shall call $A$ {\em definably amenable} if $\gen A$ carries a finitely
additive left-invariant measure $\mu$ on its definable subsets such that
$\mu(A)=1$. In this paper we shall show that in any group $G$, a definably
amenable approximate subgroup $A$  gives rise to a type-definable subgroup
$H\subseteq A^4$, such that finitely many left translates of any definable
superset of $H$ cover $A$. Hence such an approximate subgroup $A$ allows a real Lie model without expanding the language. Our proof follows
the ideas of Sanders, except that we use the measure not to define the subgroup
we obtain, but only to show that the formulas we construct in the original
language have the necessary properties. We conjecture that even without the
definable amenability assumption a suitable Lie model exists.

\smallskip

The classification of approximate subgroups of
real Lie groups is still an open problem. Since in a real Lie group any compact
neighbourhood of the identity is an approximate subgroup, in particular no
nilpotency (or even solubility) result can hold in general. We hope that under
additional model-theoretic assumptions on the original structure, a partial
classification might be easier to achieve.

\bigskip

We will end this introduction with two useful remarks. The first 
one concerns essentially the only (but crucial) use of model theory in this 
paper. The second one is an easy generalization which played a key role in the 
conclusion of \cite{BGT12}, and thus seems worth noticing. 

\smallskip

We shall assume that all structures under consideration are
\emph{$\omega^+$-saturated}, which means that any countable intersection of definable
sets is non-empty as soon as all finite subintersections are. All
non-principal ultraproducts are $\omega^+$-saturated; the compactness theorem of
model theory implies that we can replace any structure $M$ by a superstructure
$M^*$
satisfying the same first-order sentences with parameters in $M$ (an elementary
extension) which in addition is $\omega^+$-saturated.

\smallskip

As in \cite{BGT12}, the results in this paper remain true if $G$ is only a 
{\em local group}, i.e.\ a set closed under
inverse and endowed with a multiplication such that the product of up to 100
elements is well-defined and fully associative. For this, one can check 
troughout the proofs that one never needs to multiply more than 100 elements 
of $A$. 

\section{A type-definable version of Sanders' Theorem}

\defn A subset $A$ of a (local) group $G$ is {\em symmetric} if $1\in A$, and
$a^{-1}\in A$ for all $a\in A$.\par
If $K<\omega$, a symmetric subset $A$ of $G$ is a {\em $K$-approximate
subgroup} if $A^2=\{aa':a,a'\in A\}$ is contained in $K$ left cosets of
$A$. An {\em approximate subgroup} is a symmetric subset which is a
$K$-approximate subgroup for some $K<\omega$.\edefn
From a model-theoretic point of view, a definable approximate subgroup $A$
is just a symmetric {\em generic set} in $\gen A$, i.e.\ a
definable symmetric subset of $\gen A$ such that every definable subset of $\gen
A$ is covered by finitely many left translates of $A$.
\defn A definable approximate subgroup $A$ is {\em definably amenable} if there
is a left translation-invariant finitely additive measure $\mu$ on
the definable subsets of $\gen A$ with $\mu(A)=1$.\edefn
Note that by $\omega^+$-saturation, for any definable subset $X$ of
$\gen A$ there is $n<\omega$ with $X\subseteq A^n$. So if $A$ is a definably
amenable approximate subgroup, then
$\mu(X)\le\mu(A^n)<\infty$.
\bem If $\lim_{n\to\infty}\mu(A^n)<\infty$, then there is $n<\omega$ with
$A^n=\gen A$.\ebem
\bew Suppose not. Then for every $n<\omega$ there is $a_n\in A^{n+1}\setminus
A^n$. But then $(a_{3k}A:k\le n)$ is a sequence of disjoint left
translates of $A$ inside $A^{3n+2}$, whence $\mu(A^{3n+2})\ge
(n+1)\,\mu(A)=n+1$,
a contradiction.\qed

For the remainder of the paper we fix a $K$-approximate subgroup $A$ of a
(local)
group $G$, and consider the structure whose domain is $G$, with a predicate for
$A$, and with group multiplication (which is a partial map in case $G$ is only
local). We assume that $A^n\subseteq G$ for all $n<\omega$ (in fact $n\le100$
would be enough). Definability and type-definability will be with respect to
this structure.

We assume that $A$ is definably amenable with
$\lim_{n\to\infty}\mu(A^n)=\infty$. We also fix a set $E$ of size $K$ with
$A^2\subseteq EA$. 

\tats[Ruzsa's covering lemma]\label{Ruzsa}
Let $X,Y\subseteq G$ be definable such that $\mu(XY) \leq K\mu(Y)$. Then $X
\subseteq ZYY^{-1}$ for some finite $Z\subseteq X$ with $|Z| \leq K$. 
\etats
\bew If $X=\emptyset$ there is nothing to show. Otherwise, consider a
finite subset $Z\subseteq X$ such that $zY \cap z'Y = \emptyset$ for all $z
\neq z'$ in $Z$. By left invariance, 
$$|Z|\,\mu(Y)=\mu(ZY)\le\mu(XY)\le K\mu(Y),$$
so $|Z|\leq K$, and there is a maximal such $Z$. But then for any $x\in X$
there is $z\in Z$ with $zY \cap xY \neq \emptyset$ by maximality, whence $x \in
zYY^{-1}$.\qed

\defn A definable subset $B\subset\gen A$ is {\em wide in $A$} if
$A$ is covered by finitely many translates of $B$.\par
Two approximate subgroups are said to be \emph{equivalent} if each one is
contained in finitely many translates of the other.\edefn

We will sometimes make explicit the finite constants and say that $B$ is
$L$-wide in $A$, or that $A$ and $A^*$ are $L$-equivalent. 

\lem\label{wide} Let $B\subset\gen A$ be definable.\begin{enumerate}
\item If $\mu(B)>0$, then $BB^{-1}$ is wide in $A$ and symmetric.
\item If $B$ is wide in $A$ and symmetric, then $B$ is also an
approximate subgroup equivalent to $A$.\end{enumerate}\elem
\bew\begin{enumerate}
\item Clearly, $BB^{-1}$ is symmetric. Since $AB$ is a definable subset of $\gen
A$, we have $\mu(AB)<\infty$ and there is $L<\omega$ with $\mu(AB)\le L\mu(B)$.
By Fact \ref{Ruzsa}, at most $L$ translates of $BB^{-1}$ are needed to cover
$A$.
\item There is $n<\omega$ such that $B^2\subseteq A^{2n}\subseteq
E^{2n-1}A$. Suppose $Y$ is finite with $A\subseteq YB$. Then 
$$B^2\subseteq E^{2n-1}A\subseteq E^{2n-1}YB.$$
Thus $B$ is an approximate subgroup; being wide in $A$, it must be equivalent
to $A$.\qed\end{enumerate}
\defn A type-definable subgroup $H$ of a (local) group $G$ has {\em bounded}
index if there is some cardinal $\kappa$ such that in any elementary
extension the index $|G:H|$ is bounded by $\kappa$.\edefn
\bem By $\omega^+$-saturation $H$ has bounded index in $G$ if and only
if for
every definable subset $X$ of $G$ and every definable superset $Y$ of $H$,
finitely many left transates of $Y$ cover $X$.\ebem
\lem If $A$ and $A^*$ are equivalent, there exists an approximate subgroup in
which both are wide, and another one which is wide in both. In particular, $\gen
A\cap \gen {A^*}$ will have bounded index in both $\gen A$ and $\gen{A^*}$.\elem
\bew Suppose ${A^*}^2\subseteq E^*A^*$, and put $B=AA^*A$, a symmetric set
containing $A$ and $A^*$. If $A\subseteq XA^*$ and $A^*\subseteq X^*A$, then
$$B=AA^*A\subseteq XA^*A^*A\subseteq XE^*A^*A\subseteq XE^*X^*AA\subseteq
XE^*X^*EA\subseteq XE^*X^*EXA^*,$$
so $A$ and $A^*$ are wide in $B$. Moreover
$$B^2\subseteq XE^*X^*EAB=XE^*X^*EAAA^*A\subseteq XE^*X^*E^2AA^*A=XE^*X^*E^2B,$$
so $B$ is also an approximate subgroup. As $\gen A$ and $\gen{A^*}$ have
bounded index in $\gen B$, the intersection $\gen A\cap\gen{A^*}$ has bounded
index in $\gen B$, and thus in $\gen A$ and in $\gen A^*$.

Now note that $\gen A\cap\gen{A^*}=\bigcup_{n<\omega}(A^n\cap{A^*}^n)$. By
$\omega^+$-saturation there is $n<\omega$ such that finitely many
translates of
$A^n\cap{A^*}^n$ cover $B$. So $A^n\cap{A^*}^n$ is wide in $B$, whence in $A$
and in $A^*$.\qed

We now turn to the main result. We shall need the following Lemma due to
Sanders.
\lem\label{lmm_f(t)} Let $f :\,]0,1]\to[1,K]$ and $\epsilon > 0$. Then there
exists $n<\omega$ depending only on $K, \epsilon$ and $t > 1 / (2K)^{2^n -1}$
such that 
$$f(\frac{t^2}{2K} )\geq(1-\epsilon)f(t).$$\elem
\bew Define a sequence $(t_n)$ by $t_0 = 1$ and $t_{k+1} =
\frac{t_k^2}{2K}$, so $t_n = 1 / (2K)^{2^n-1}$. 

For $n<\omega$ suppose that for all $i <n$ we have $f(t_{i+1}) <
(1-\epsilon) f(t_i)$. Then 
$$f(t_n) < (1-\epsilon)^n f(t_0) \leq (1-\epsilon)^n K.$$
But $f(t_n) \geq 1$, so if $n<\omega$ is such that
$(1-\epsilon)^n K < 1$, there must be some $i < n$ with $f(t_{i+1})
\geq (1-\epsilon) f(t_i)$.\qed

\satz\label{thm_Sanders} Let $A$ be a $K$-approximate subgroup. For any
$m<\omega$ there is a definable $L$-wide approximate subgroup $S$ with
$S^m\subseteq A^4$, where $L$ depends only on $K$ and $m$.\esatz
\bew Let us show first that if $B\subseteq A$ is definable with $\mu(B)\ge
t\mu(A)$ for some $0<t\le 1$ and $s=\frac t{2K}$, then $A$ is covered by
$N=\llcorner\frac1s\lrcorner$ translates of
$$X=\{g\in A^2:\mu(gB\cap B)\ge st\mu(A)\}$$
by elements of $A$. So suppose not. Then inductively we find a
sequence $(g_i:i\le N)$ of elements of $A$ such that $\mu(g_iB\cap
g_jB)<st\mu(A)$ for all $i<j\le N$, since for $j\le N$ the set
$\bigcup_{i<j}g_iX$ cannot cover $A$. But then
$$\begin{aligned}K\mu(A)&\ge\mu(A^2)\ge\mu(\bigcup_{i\le N}g_iB)\\
&\ge(N+1)\mu(B)-\sum_{i<j\le N}\mu(g_iB\cap g_jB)\\
&>(N+1)t\mu(A)-\frac{N(N+1)}2 st\mu(A)=(1-N\frac s2)(N+1)t\mu(A)\\
&\ge(1-\frac1s\frac s2)\frac1s
t\mu(A)=\frac12\frac{2K}tt\mu(A)=K\mu(A),\end{aligned}$$
a contradiction. 

\smallskip

However, as $\mu$ is not supposed to be definable, $X$ need not be definable
either. We shall hence look for definable sets with similar properties. 
To this end, consider the following conditions $P_n^t(X)$ on definable
subsets of $A$, for $n<\omega$ and $0<t\le 1$:
\begin{itemize}
\item $P_0^t(B)$ if $B\not=\emptyset$.
\item $P_{n+1}^t(B)$ if $P_n^t(B)$, and $A$ is
covered by $\llcorner\frac{2K}t\lrcorner$
translates of 
$$X_{n+1}^t(B)=\{g\in A^2:P_n^{t^2/2K}(gB\cap B) \ \textrm{and} \
P_n^{t^2/2K}(g^{-1}B\cap B)\}.$$\end{itemize}
Clearly, if $(B_x)_x$ is a family of uniformly definable subsets of
$A$, then $P_n^t(B_x)$ is definable by a formula $\theta_n^t(x)$ for all
$n<\omega$ and $0<t\le1$. As $X_{n+1}^t(B)\subseteq A^2$, the translating
elements for the covering of $A$ must come from $A^3$, so the $P_n^t$ are
definable even in a local group (where we can only quantify over
finite powers of $A$).

For $0<t\le1$ we consider the family $\BB_t$ of definable subsets $B$ of $A$
with $P_n^t(B)$ for all $n<\omega$. The first paragraph implies
inductively that for definable $B\subseteq A$, if $\mu(B)\ge t\mu(A)$
then $P_n^t(B)$ holds, whence $B\in\BB_t$. In particular, $A \in \BB_t$ so
$\BB_t$ is non-empty. Note that $P_n^t$ implies $P_n^{t'}$ for $t\ge t'$, so
$\BB_t\subseteq\BB_{t'}$.

\smallskip

Define a function $f:\,]0,1]\to\mathbb R$ by 
$$f(t)=\inf\{\frac{\mu(BA)}{\mu(A)}:B\in\BB_t\}.$$
Fix $\epsilon > 0$. Since $1\le f(t)\le K$ for all $0<t\le1$, by
Lemma~\ref{lmm_f(t)} there is $t>0$ depending only on $K$ and $\epsilon$
such that 
$$f(\frac{t^2}{2K})\ge(1-\epsilon)f(t).$$
Choose $B\in\BB_t$ with $\frac{\mu(BA)}{\mu(A)}\le(1+\epsilon)f(t)$. Put 
$$X_n=X_n^t(B)=\{g\in A^2:P_n^{t^2/2K}(gB\cap B)\ \textrm{and} \
P_n^{t^2/2K}(g^{-1}B\cap B)\}$$
and $X=\bigcap_{n<\omega}X_n$. Then $X_n$ is symmetric, $X_{n+1}\subseteq X_n$
and $\llcorner\frac{2K}t\lrcorner$ translates of $X_n$ cover
$A$, for all $n<\omega$.
By $\omega^+$-saturation, $\llcorner\frac{2K}t\lrcorner$ translates of
$X$ cover $A$, so
$X$ is nonempty. Moreover, for $g\in X$ we have $gB\cap
B\in\BB_{t^2/2K}$, whence 
$$\begin{aligned}\mu(gBA\cap BA)&\ge\mu((gB\cap B)A)\ge
f(\frac{t^2}{2K})\mu(A)\\
&\ge(1-\epsilon)f(t)\mu(A)\ge\frac{1-\epsilon}{1+\epsilon}\,\mu(BA).
\end{aligned}$$
Hence for $g\in X$,
$$\mu(gBA\triangle
BA)\le\frac{4\,\epsilon}{1+\epsilon}\,\mu(BA)<4\,\epsilon\,\mu(BA).$$
It follows that for $g_1,\ldots,g_m\in X$,
$$\begin{aligned}\mu(g_1\cdots g_mBA&\triangle BA)\\
&\le\mu((BA\triangle g_1BA)\cup g_1(BA\triangle g_2BA)\cup\cdots\cup g_1\cdots
g_{m-1}(BA\triangle g_mBA))\\
&\le\mu(BA\triangle g_1BA)+\mu(BA\triangle g_2BA)+\cdots+\mu(BA\triangle
g_mBA)\\
&<4\,m\,\epsilon\,\mu(BA).\end{aligned}$$
In particular, if $\epsilon\le\frac1{4m}$, then $g_1\cdots g_m BA\cap
BA\not=\emptyset$, whence $X^m\subseteq A^4$.
By $\omega^+$-saturation there is $n<\omega$ such that $X_n^m\subseteq
A^4$. Note that
$S := X_n$ is $\llcorner\frac{2K}t\lrcorner$-wide in $A$, and thus an
approximate subgroup equivalent to $A$ by Lemma~\ref{wide}.\qed

\kor\label{thm_full} There is a type-definable subgroup $H\subseteq A^4$
such that every definable superset of $H$ contained in $\gen A$ is wide in
$A$.\ekor
\bew Put $S_0=A$ and apply inductively Theorem \ref{thm_Sanders} for $m=8$ with
$S_i$ instead of $A$, in order to obtain a sequence of approximate
subgroups $(S_i:i<\omega)$ with $S_{i+1}$ wide in $S_i$ (whence in $A$) and
$S_{i+1}^8\subseteq S_i^4$. Then $H=\bigcap_{i<\omega}S_i^4$ is a type-definable
subgroup of $A^4$. Any definable superset of $H$ must contain some $S_i^4$ by
$\omega^+$-saturation, and hence be wide in $A$.\qed

\section{Normality}
Since we want to consider the quotient $\gen A/H$, we shall look for
a stronger version of Theorem \ref{thm_Sanders} where $H$ will be normal.

\lem\label{lmm_normal}
Let $X_1,..,X_n$ be definable subsets of $A$ with $N_i\,\mu(X_i) \geq
\mu(A)$ for some $N_i<\omega$. Then there is a definable $D \subseteq
A$ such that 
\[D^{-1}D \subseteq (X_1X_1^{-1})^2 \cap ... \cap (X_nX_n^{-1})^2
\quad \textrm{and} \quad K^{n-1}N_1...N_n\,\mu(D) \geq \mu(A).\]
\elem
\bew 

Since $\mu(AX_2) \leq K\,\mu(A) \leq KN_2\,\mu(X_2)$, by
Fact~\ref{Ruzsa} there are $g_1,\ldots,g_{KN_2}$ such that
\[A \subseteq \bigcup_{i=1}^{KN_2} g_i X_2X_2^{-1}.\]  
Then there is an $i$ such that
$$KN_1N_2\,\mu(X_1 \cap g_iX_2X_2^{-1}) \geq \mu(A).$$
We set $D_0= X_1 \cap
g_iX_2X_2^{-1}$ and note that $D_0^{-1}D_0 \subseteq X_1^{-1}X_1 \cap
(X_2X_2^{-1})^2$.

\smallskip

Then we can iterate the construction, replacing $X_1$ by $D_0$
and $X_2$ by $X_3$. Inductively we obtain a suitable $D$
with $K^{n-1}N_1...N_n\,\mu(D) \geq \mu(A)$ such that 
$$D^{-1}D \subseteq X_1^{-1}X_1 \cap (X_2X_2^{-1})^2 \cap ... \cap
(X_nX_n^{-1})^2.$$
Notice that $D^{-1}D$ is $K^{n-1}N_1...N_n$-wide in $A$ by Fact~\ref{Ruzsa}. 
\qed

\satz\label{thm_normal} Let $A$ be a $K$-approximate subgroup, and $R$ a
definable
$N$-wide symmetric subset with $R^4 \subseteq A^4$. Then there exists a
definable
$L$-wide symmetric subset $S$ with $(S^8)^A \subseteq R^4$,
where $L$ depends only on $K$ and $N$.\esatz
\bew If $A\subseteq XR$, then 
$$R^2\subseteq A^4\subseteq E^3A\subseteq E^3XR,$$
so $R$ is a $K^3N$-approximate subgroup. Theorem
\ref{thm_Sanders} yields the existence of some $T\subseteq R^4$
equivalent to $R$ with $T^{48}\subseteq R^4$. Then $T$ is wide in $A$ and there
exists $n<\omega$ depending only on $K$ and $N$ and some elements $a_i$
of $A$ such that 
\[A \subseteq \bigcup_{i=1}^n a_iT.\]
Consider the measure $\bar{\mu}$ on definable subsets of $\gen A$ defined by
\[\bar{\mu}(X) := \frac{1}{n} \sum_{i=1}^n \mu(Xa_i) .\] 
Clearly $\bar\mu$ is left translation invariant, we have
$$\bar{\mu}(A) = \frac{1}{n} \sum_{i=1}^n \mu(Aa_i)\le\frac nn\,\mu(A^2)\leq K
\mu(A),$$
and 
$$\bar{\mu}(a_iTa_i^{-1}) \geq\frac{1}{n}\mu(T) \geq\frac{1}{n^2}
\mu(A) \geq\frac{1}{Kn^2}\bar{\mu}(A).$$ 
Since all the $a_iTa_i^{-1}$ are subsets of $A^6$ and
$$K^6n^2\bar{\mu}(a_iTa_i^{-1})\geq K^5\bar\mu(A)\ge \bar{\mu}(A^6),$$
Lemma~\ref{lmm_normal} applied to the $K^6$-approximate subgroup $A^6$ yields a
subset $D\subseteq A^6$ with 
$$(K^6)^{n-1}(K^6n^2)^n\bar\mu(D)\ge\bar\mu(A^6)$$
such that for $i=1,2,\ldots,n$ we have
\[S := D^{-1}D\subseteq a_iT^4a_i^{-1}.\]
Then $S$ is symmetric, wide in $A$ and $S^{a_i}\subseteq T^4$ for
$i=1,\ldots,n$.
Since $A\subseteq \bigcup a_iT$, this means that $S^A \subseteq T^6$, so
$(S^8)^A\subseteq T^{48}\subseteq R^4$.\qed

\kor There is a type-definable normal subgroup $H$ of $\gen A$
contained in $A^4$ such that every definable superset of $H$ contained in $\gen
A$ is wide in $A$.\ekor
\bew As Corollary \ref{thm_full}, using Theorem \ref{thm_normal} instead of
Theorem \ref{thm_Sanders}.\qed

\end{document}